\documentclass[reprint,twocolumn,amsmath,amssymb,aps]{revtex4-1}
\pdfoutput=1 
\usepackage{graphicx}
\usepackage{dcolumn}
\usepackage{bm}
\usepackage{hyperref}

\usepackage{xcolor}

\DeclareMathOperator*{\argmin}{arg\,min}

\begin{document}

\title{Probabilistic Behavioral Distance and Tuning\\ - Reducing and aggregating complex systems}

\author{Frank Hellmann}
\email{hellmann@pik-potsdam.de}
\affiliation{Potsdam Institute for Climate Impact Research}

\author{Ekaterina Zolotarevskaia}
\email{zolotarevskaia@pik-potsdam.de}
\affiliation{Potsdam Institute for Climate Impact Research}

\author{Jürgen Kurths}
\affiliation{Potsdam Institute for Climate Impact Research}

\author{Jörg Raisch}
\affiliation{TU Berlin}

\date{\today}

\begin{abstract}
Given a complex system with a given interface to the rest of the world, what does it mean for a the system to behave close to a simpler specification describing the behavior at the interface? We give several definitions for useful notions of distances between a complex system and a specification by combining a behavioral and probabilistic perspective. These distances can be used to tune a complex system to a specification. We show that our approach can successfully tune non-linear networked systems to behave like much smaller networks, allowing us to aggregate large sub-networks into one or two effective nodes. Finally, we discuss similarities and differences between our approach and $H_\infty$ model reduction.
\end{abstract}

\maketitle

\section{Introduction}

The behavioral approach to the theory of dynamical systems focuses on the possible trajectories that a system can exhibit \cite{willems1989models, willems1997introduction}. In the context of control theory these could be for example the inputs and outputs of a closed loop system. This perspective allows one to speak about specification and system behavior on the same level. The basic notion is: A dynamical system is a set of possible trajectories, a specification is a set of permitted trajectories and a dynamical system satisfies a specification if every possible trajectory is permissible.

This paper introduces several notions of distance between a system and a specification. These distance concepts are probabilistic, meaning they can be usefully approximated. Obtaining meaningful and usable notions of the distance of system and specification enables us to tune complex systems to satisfy specifications at least approximately. The context we have in mind is tuning the control of a subsystem of a larger network to present a unified aggregate behavior towards the "rest of the world". In this setting the links towards the larger network act as inputs/outputs. Another application would be the hierarchical decomposition of a larger control task, where the specification of lower levels in the hierarchy serve as systems for the higher level \cite{moor2003admissibility}.

Two different distance notions are introduced, one suited for optimizing systems towards a specification, and one for probabilisitically validating their compliance. We show that standard tools from differential programming and non-linear optimization can be used successfully to tune complex systems by minimizing the sampled approximation of the distance by tuning several complex non-linear dynamical networks.

The method developed here is designed with applications to power grids in mind. In parallel to this paper we present a software stack based on the capabilities of the Julia language \cite{rackauckas2017differentialequations, rackauckas2019diffeqflux}, that allows tuning dynamical properties in power systems \cite{buttner2021stack, lindner2021networkdynamics, plietzsch2021powerdynamics}. In this context a behavioral approach is particularly natural. Typically the regulatory specification of how the system should behave is not given in terms of a precise set of differential equations, but rather by describing general properties of the trajectories. For example the system limits for the rate of the change of frequency (RoCoF) at which disconnections of major generators will occur, and thus cascading blackouts become likely, is defined by the European transmission system operators by the conditions that the moving average of the RoCoF stays within: $\pm 2$Hz/s for a 500ms window, $\pm 1,5$Hz/s for  a 1000ms window and $\pm 1,25$Hz/s for a 2000ms window\cite{entsoe}. Any frequency trajectory that stays inside this curve is considered acceptable. At the same time, the demands, perturbations and faults that the power grid experiences are varied and random, necessitating a probabilistic approach. This has long been standard in the analysis of static properties, e.g. \cite{borkowska1974probabilistic, anders1989probability}, but is increasingly also used to systematically understand dynamic aspects of the system \cite{menck2014dead, hellmann2016survivability, hellmann2020network, liemann2020probabilistic}.



\section{Systems and specifications}

The highly abstract definitions of behavioral dynamical systems theory are hard to work with directly. The setting of this paper is to consider behaviors given by parametrized input-output differential equations. Fix some time interval $T = [0, t_{\text{final}}]$. Then denote the trajectories of the internal states $x \in X^T$, the input states $i \in I^T$ , and the output states $o \in O^T$. \textit{We will always use $x$, $i$ and $o$ to refer to the function, and $x(t) \in X$ to refer to a concrete value.} The dynamical system is then given by specifying the dynamics $f$, the output function $g$, and the initial conditions $x_0$, all of which can depend on parameters $p \in \mathcal{P}$. Finally the inputs are restricted to some set $\mathcal{B}_i$, the equations then are:

\begin{align}\label{eq:IO-ODE}
    \dot x(t) &= f(x(t),i(t),p)\\
    o(t) &= g(x(t),p) \nonumber\\
     x(0) &= x_0(p) \nonumber\\
  p \in \mathcal{P} &\text{ and }
   i \in \mathcal{B}_i \subset I^T\nonumber
\end{align}

We always assume that these equations can be integrated for the time period $T$. The set of possible trajectories of this set of equations is parametrized by $\mathcal{P} \times \mathcal{B}_i$.

In what follows it will be important to distinguish between a system with parameters, a system without parameters and a specification. Even though all behaviors we consider will be of the IO form \eqref{eq:IO-ODE}, we will denote specifications using $z \in Z^T$ for internal state, and $q$ for parameters, and we will denote the complete system of the specification as $\mathcal{C}$:

\begin{align}\label{eq:spec}
    \mathcal{C} : \dot z(t) &= f^\mathcal{C}(z(t),i(t),q) \\
    o^\mathcal{C}(t) &= g^\mathcal{C}(z(t),q) \nonumber\\
    z(0) &= z_0(p) \nonumber\\ \text{ and } & \nonumber q \in \mathcal{Q}\nonumber 
\end{align}

A parametrized system is denoted by $\mathcal{S}$:

\begin{align}\label{eq:system}
    \mathcal{S} : \dot x(t) &= f^\mathcal{S}(x(t),i(t),p)\\
    o^\mathcal{S}(t) &= g^\mathcal{S}(x(t),p) \nonumber\\
     x(0) &= x_0(p) \nonumber\\ \text{ and } & \nonumber p \in \mathcal{P}\nonumber
\end{align}

Finally, given a system with no parameter freedom (e.g. $|\mathcal{P}| = 1$), we call the system unparametrized and denote it as $\overline{\mathcal{S}}$. The behavior of the system is then determined entirely by the set of inputs $\mathcal{B}_i$. We denote the unparametrized system obtained by fixing the parameter of a system $\mathcal{S}$ at some $p \in \mathcal{P}$ by $\mathcal{S}|p$.

\begin{figure}[h!b]
\includegraphics[width=0.49\columnwidth]{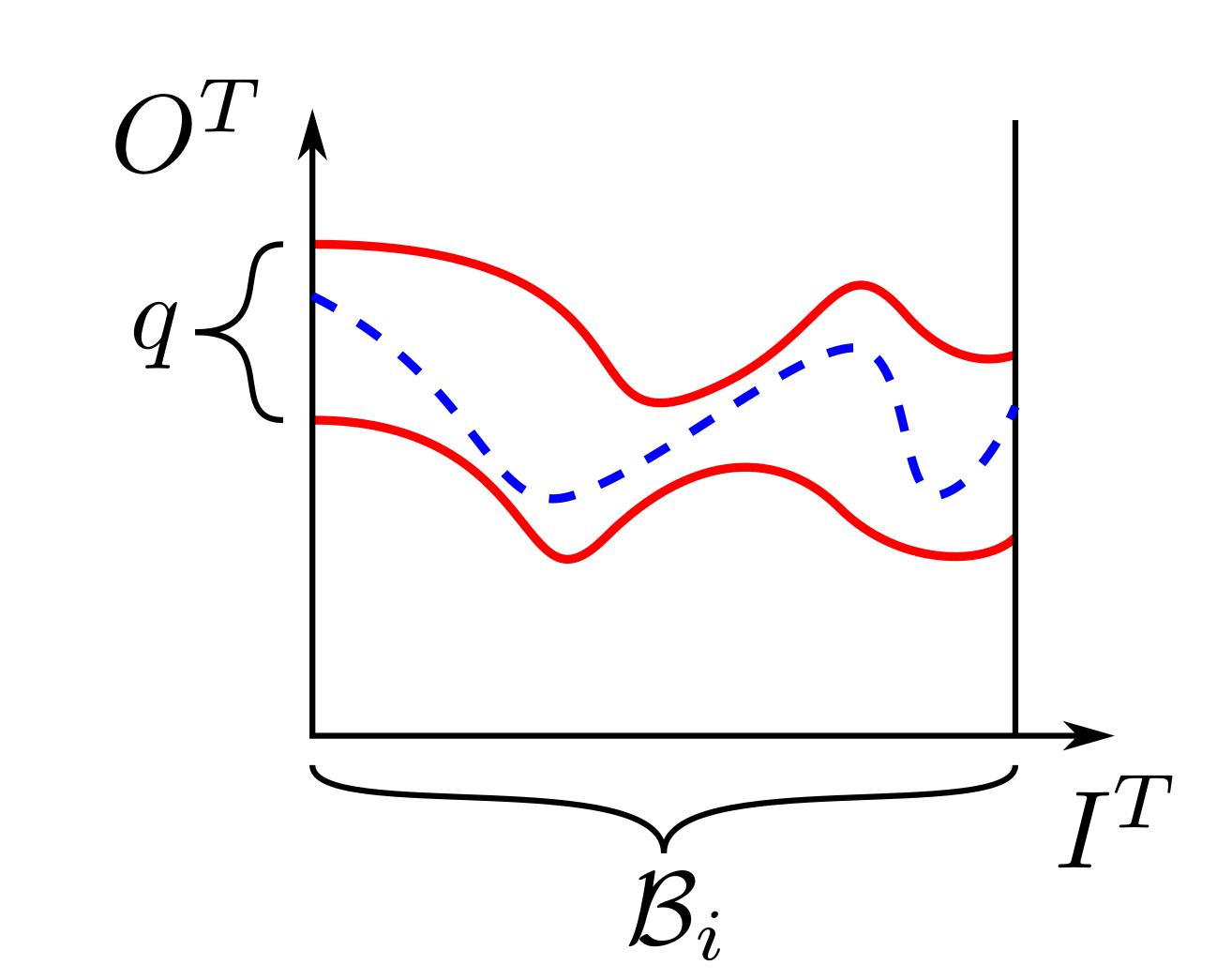}
\includegraphics[width=0.49\columnwidth]{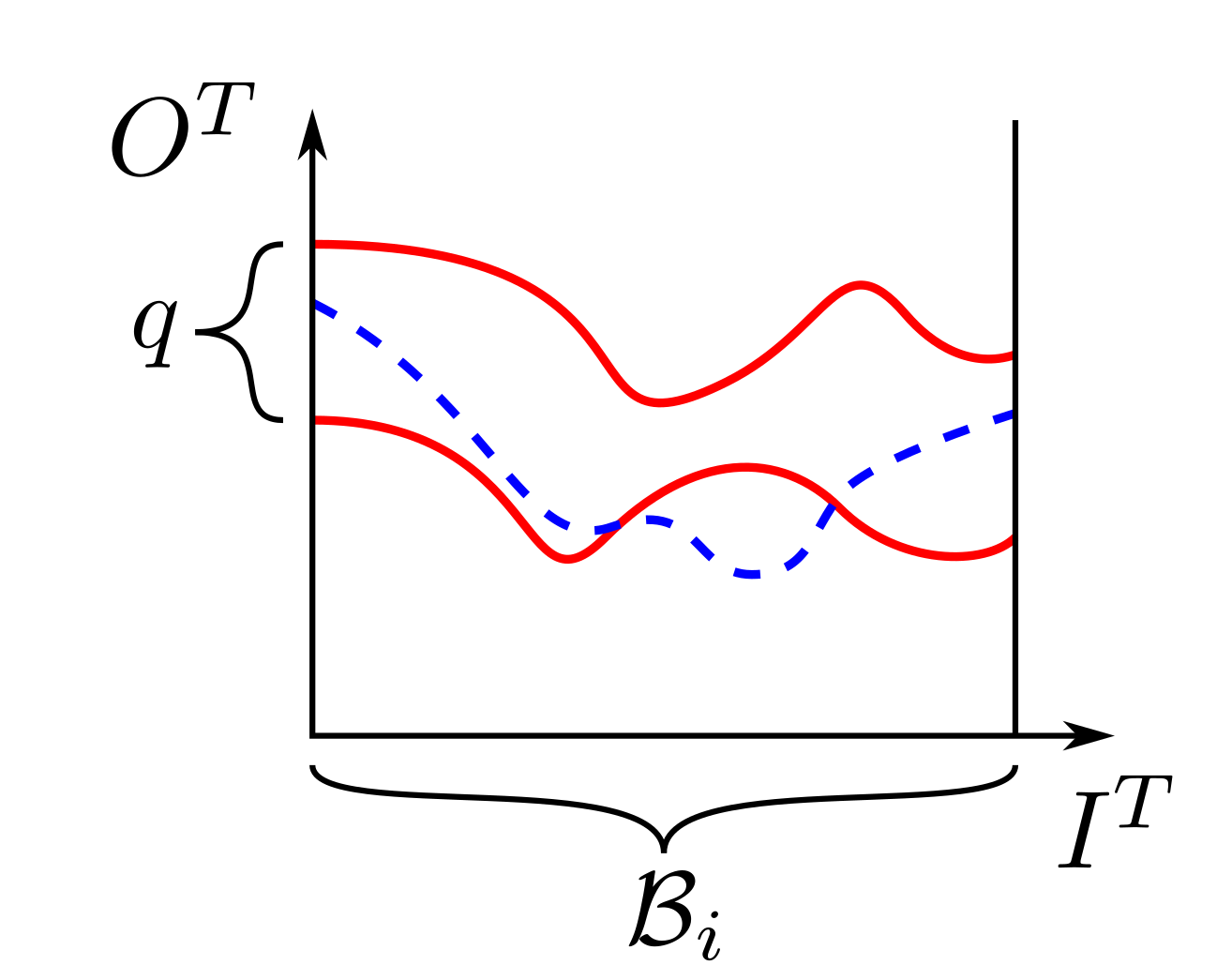}
\caption{On the horizontal axis we have the set of functions $I^T$, on the vertical axis the set of functions $O^T$. The solid red lines bound the range of outputs $o$ that can occur in the specification given a specific $i$. This range is parametrized by $q$ in equation \eqref{eq:spec}. An input output system defines a graph, as it maps a given input to an output. On the left, the system with the blue dashed graph satisfies the specification as, for every $i$ there is a $q$ such that the $o$ of the specification matches the $o$ of the system. On the right there are some inputs for which no such $q$ exists, and the system does not satisfy the specification. \label{fig:satisfy spec}}
\end{figure}

Using these notions of behavior, we can now see explicitly what it means in our case for a system to satisfy a specification. Take a specification $\mathcal{C}$ parametrized by $\mathcal{Q}$. Take an unparametrized system $\overline{\mathcal{S}}$. The output $o$ is completely determined by the input $i$. $\overline{\mathcal{S}}$ satisfies the specification if for every input $i$ there is a $q[i] \in \mathcal{Q}$, such that the specification system $\mathcal{C}$ will produce the output $o$. This is illustrated in Figure~\ref{fig:satisfy spec}.

\textbf{Remark:} \textit{Note that it is not necessary that there is one set of parameters $q$ that matches the system behavior for all inputs for $\overline{\mathcal{S}}$ to satisfy the specification. In particular we do not require that the two ODEs transfer functions match. The behavioral condition is strictly weaker.}

Given a system with parameters, the tuning problem we want to address is to find a $p$ such that $\mathcal{S}|_p$ satisfies a specification. To do so we will now introduce several notions of distance between system and specification.

\section{Distance to the specification}

In practice it might often be impossible but also unnecessary to satisfy the specification exactly. This can mean both, that it is acceptable to fail for some inputs, or that the outputs are not exactly but only approximately the same. Our goal will be to get the system to be close to the specification. To formalize this idea we introduce a notion of the distance of a system to the specification.

The first important ingredient for this is a distance on the set of output functions, $\Delta(o_1, o_2)$. In what follows we will always take the square of the $L_2$-norm:

\begin{align}
\Delta(o_1, o_2) = \int_T \|o_1(t) - o_2(t)\|^2 \text{d}t
\end{align}

As noted above, the outputs are functions of the inputs and the parameters. We will write $o[i, p]$ or $o[i]$ if no parameters are present. At fixed input $i$ and parameter $p$, $o[i, p]$ is a function of time. Now given an unparametrized system $\overline{\mathcal{S}}$ and a specification $\mathcal{C}$ of the form \eqref{eq:system} and \eqref{eq:spec} and a given input $i$ we can consider

\begin{align}
\min_q \Delta(o^{\overline{\mathcal{S}}}[i], o^{\mathcal{C}}[i, q]) \label{eq:distance_at_i}
\end{align}

as a distance of the system to the specification at input $i$.

\begin{figure}[hb]
\includegraphics[width=0.69\columnwidth]{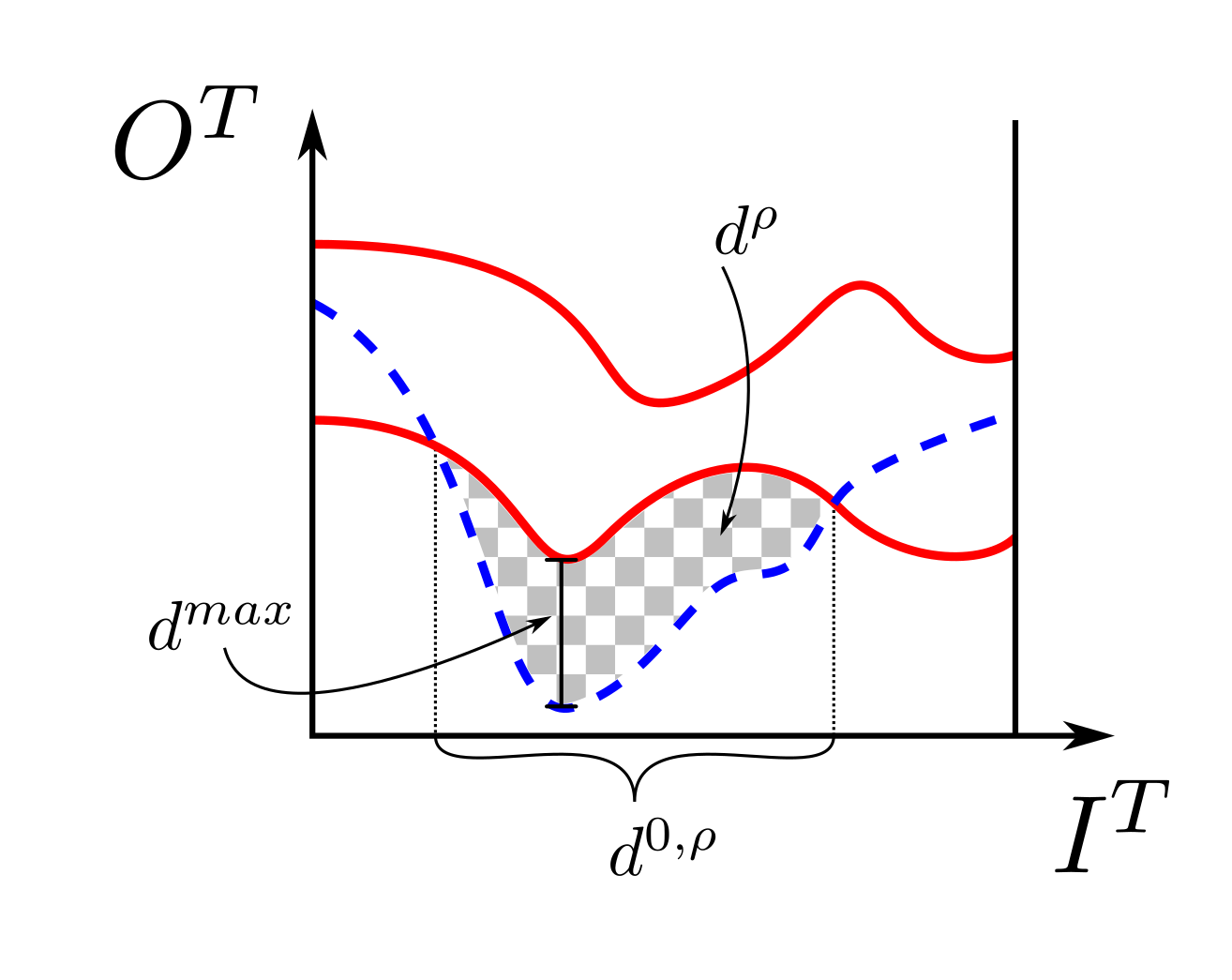}
 \caption{Using the same way to illustrate a specification and a system as in Figure~\ref{fig:satisfy spec} we illustrate the three distance measures. The area in between the closest elements of the specification and the system output is checkered, and provides the distance $d^\rho$. The size of the range of $i$ for which the specification is violated is drawn on the horizontal axis and corresponds to $d^{0,\rho}$, that is $d^{\epsilon,\rho}$ for $\epsilon = 0$, finally the largest distance between the specification and the system corresponds to $d^{max}$. \label{fig:distances}}
\end{figure}

In order to understand how much the outputs typically diverge we need to provide information on what inputs the system typically encounters. This can be formalized by providing a probability measure $\rho$ on $\mathcal{B}_i$. The main distance measure we will investigate in what follows is then given by taking the average distance of the system to the specification in the sense of \eqref{eq:distance_at_i}:

\begin{align}\label{eq:distance_rho}
d^{\rho}(\overline{\mathcal{S}}, \mathcal{C}) &= \int_{\mathcal{B}_i} \min_q \Delta(o^{\overline{\mathcal{S}}}[i], o^\mathcal{C}[i, q]) \rho\\
&= E_\rho\left [ \min_q \Delta(o^{\overline{\mathcal{S}}}[i], o^\mathcal{C}[i, q]) \right ]\nonumber
\end{align}

where $E_\rho$ denotes the expectation value.

Note that $d^{\rho}(\overline{\mathcal{S}}, \mathcal{C}) = 0$ implies that the system satisfies the spec for all $i$ up to a set of measure zero in $\rho$. Importantly, if the expectation is small we also know that the output of the system is close to an output allowable according to the specification for most inputs.

A second complementary distance measure that will be useful for validating the performance of a system in a more rigorous way is the fraction of inputs for which the distance of the system to the specification is larger than a margin $\epsilon$. Let $\Theta_\epsilon:\mathbb{R}^n\rightarrow\{0,1\}$ be the shifted Heaviside step function with $\Theta_\epsilon(z) = 1$ if $z > \epsilon$ and $0$ otherwise. Then we have:

\begin{align}
d^{\rho, \epsilon}(\overline{\mathcal{S}}, \mathcal{C}) &= \int_{\mathcal{B}_i} \Theta_\epsilon( \min_q \Delta(o^{\overline{\mathcal{S}}}[i], o^\mathcal{C}[i, q]) ) \rho\nonumber\\
&= E_\rho \left [ \Theta_\epsilon ( \min_q \Delta(o^{\overline{\mathcal{S}}}[i], o^\mathcal{C}[i, q] )) \right ] \label{eq:distance_epsilon}
\end{align}

As $\rho$ is a probability measure, $d^{\rho, \epsilon}$ will vary between $0$ and $1$. Note that $d^{\rho, \epsilon} = 0$ does not guarantee that the system satisfies the specification. It only implies that almost all inputs produce an output that differs from the specification by less than or equal to $\epsilon$. Only for $\epsilon = 0$ do we guarantee that the specification is exactly satisfied for almost all $i$. Note that these two distances are genuinely complementary and do not coincide.

Finally we note that there is a natural distance that does not depend on a probability measure, but just on some norm on $\mathcal{B}_i$, namely the maximum distance to the specification given a fixed norm of the input:

\begin{align} \label{eq:distance_max}
d^{\max}(\overline{\mathcal{S}}, \mathcal{C}) = \max_{i : \|i\| = 1} \min_q \Delta(o^{\overline{\mathcal{S}}}[i], o^\mathcal{C}[i, q])
\end{align}

This notion of distance closely resembles an operator norm. Whereas $d^\rho$ and $d^{\rho, \epsilon}$ capture the typical performance of the system, $d^{\max}$ is concerned with the worst case performance only. We will return to this in Section~\ref{sec:h-inf} where $d^{\max}$ will allow us to connect and contrast the perspective taken here to conventional notions of the $H_\infty$ operator norm and model reduction.

\section{Sampling based approximations}

The distances introduced above require a probability distribution on a set of inputs. Probability distributions on spaces of trajectories in time are called stochastic processes. The distances we introduced are thus given by expectation values of stochastic processes, and that in turn means they can be approximated by evaluating them on a sample $B_i \subset \mathcal{B}_i$, a set of realizations of the process.

Note that, while the various kinds of stochastic differential equations are the most familiar and best studied class of stochastic processes, they are not necessarily the best suited for the differential equation setting we study here. Instead we will make use of smooth random functions or random ODEs with smooth solutions \cite{filip2019smooth}. We will give an example of this below.

Now given a sample $B_i$ of realizations of the stochastic process $\rho$, with cardinality $|B_i|$, we can introduce the estimators

\begin{align}
\hat d^{\rho} &= \frac1{|B_i|}\sum_{i \in B_i} \min_p \Delta(o^{\overline{\mathcal{S}}}[i], o^\mathcal{C}[i, q])\label{eq:estimator_distance_rho}
\end{align}

for \eqref{eq:distance_rho} and

\begin{align}
\hat d^{\rho, \epsilon} &= \frac1{|B_i|}\sum_{i \in B_i} \Theta_\epsilon ( \min_q \Delta(o^{\overline{\mathcal{S}}}[i], o^\mathcal{C}[i, q] ))\label{eq:estimator_distance_epsilon}
\end{align}

for \eqref{eq:distance_epsilon}. It is challenging to understand the quality of the first estimator rigorously without further information on the summand. However, the summand for the second estimator is either $0$ or $1$. This means we can interpret this sample as a Bernoulli trial and we can use the standard center point corrections and confidence intervals \cite{agresti1998approximate}. For example, the "add two successes and failures" 95\% confidence interval is given in terms of $\tilde{n} = |B_i| + 4$ and $$\tilde d = \frac{|B_i| \hat d^{\rho, \epsilon} + 2} {\tilde{n}}$$ by:

\begin{align}
d^{\rho, \epsilon} \approx \tilde d \pm 2 \sqrt{\frac{\tilde d(1 - \tilde d)}{\tilde n}}.
\end{align}

In practice it will typically not be possible to evaluate the minima in \eqref{eq:estimator_distance_epsilon} exactly. Note however, that this leads to an overestimation of the distance. Thus in practice one always obtains somewhat conservative estimators that overestimate the distance of the system to the specification. The confidence interval is most accurately thought of as the confidence interval for an estimator for an upper bound on the distance $d^{\rho, \epsilon}$.

While it is easier to provide rigorous statements on the quality of the estimator for $d^{\rho, \epsilon}$, the estimate for $d^{\rho}$ has the advantage of not requiring an appropriate choice of $\epsilon$ and of varying smoothly as the fit between system and specification changes. This makes it better suited as a basis for tuning the system using optimization techniques.

\section{Tuning the system}

Given a parametrized system of the form \eqref{eq:system}, we can formulate an optimization problem to find the set of parameters for which the system is the closest to the specification. Recall that $d^{\rho}$, as defined in \eqref{eq:distance_rho}, is the expectation value of a minimum. Optimizing it is thus a non-linear 2-stage stochastic programming problem. By using the approximation $\hat d^{\rho}$ we can explicitly give a discretized extensive form.

Recall that we denote as $\mathcal S|{p}$ the unparametrized system obtained by setting the parameters of $\mathcal{S}$ to $p$. Thus we want to find:

\begin{align}
p_{tuned} = \argmin_{p} \hat d^{\rho}(\mathcal S|{p}, \mathcal{C})
\end{align}

In the expectation value underlying $d^{\rho}$ there is a minimization for each input. To make this explicit, we will denote the parameter $q$ of the specification $\mathcal{C}$ for a given input $i$ as $q_i$. Then we can exchange the order of the sum and minimization:

\begin{align} \label{eq:joint-optim}
& p_{tuned} \nonumber\\
&= \argmin_{p} \frac1{|B_i|}\sum_{i \in B_i} \min_{q_i} \Delta(o^\mathcal{S}[i, p], o^\mathcal{C}[i, q_i]) \nonumber\\
&= \argmin_{p} \min_{\{q_i\}} \frac1{|B_i|} \sum_{i \in B_i} \Delta(o^\mathcal{S}[i,p], o^\mathcal{C}[i, q_i]).
\end{align}

This is as a joint parameter optimization in $p$ and the set $\{q_i\}$ of a large differential equation with a trajectory based target function. This type of optimizations can be implemented in a straightforward manner using DiffEqFlux \cite{rackauckas2019diffeqflux} in the Julia language \cite{bezanson2017julia}. Our implementation is available at \url{https://github.com/PIK-ICoNe/ProBeTune.jl}. Thanks to the ability to differentiate through ODE solvers a wide range of optimizers are available to perform this optimization. Crucially, the distance measure is designed in such a way that we can perform a joint optimization rather than having to perform an optimization of an optimum, a much harder problem.

\section{Demonstration for a non-linear network}

We will demonstrate the distances and their tuning, by considering two paradigmatic examples of complex non-linear dynamical networks connected at one node to the outside world. We will tune them to react to outside inputs like a specification given by a vastly simpler network, thus demonstrating that probabilistic behavioral tuning can aggregate complex networks. The two systems we will consider are a diffusively coupled network with tunable non-linear stabilizing forces and a system of Kuramoto oscillators with inertia.

\subsection{Diffusive non-linear network}

Consider the networked system $\mathcal{A}^{N}$ with $N$ nodes denoted by $n = 1 \dots N$:

\begin{gather}\label{eq:networked system}
\dot x_n = -x_n - p_n x_n^3 + \sum_{m = 1}^N A_{nm} (x_n - x_m) + \delta_{n1} (x_n - i),\nonumber\\
o(t) = i(t) - x_1(t)\\ \nonumber 
x_n(0) = 0 \nonumber\\
p_n \in \mathbb R^+ \nonumber
\end{gather}

\begin{gather}
\mathcal{B}_i = \left\{i : i(t) = \text{Re} \sum_{l = 0}^L a_l e^{i (2\pi l t + \theta_l)}\right\}
\end{gather}

for some fixed graph with adjacency matrix $A_{nm}$. The trajectories of this system are always bounded and therefore smooth solutions always exist. We can specify a probability distribution on $\mathcal{B}_i$ by specifying a probability for the parameters $a_l \in \mathbb R$ and $\theta_l \in [0, 2\pi)$.

Our goal now is to tune a full system with $N$ nodes to behave, as far as the input-output relationship is concerned, like a two node system $\mathcal{A}^2$. In the concrete example we choose the Barabasi-Albert model \cite{barabasi1999emergence} with $N =10$ to generate $A_{nm}$, as such a scale-free network has a rich irregular structure.

Thus we have $\mathcal{S} = \mathcal{A}^{10}$ and $\mathcal{C} = \mathcal{A}^{2}$. In \eqref{eq:joint-optim} this implies that we will jointly optimize over $({\mathbb R^+})^{10 + 2|B_i|}$. That is, the system $\mathcal{A}^{N}$ with $N$ parameters and one copy of the system $\mathcal{A}^{2}$ for each element of the sample $B_i \subset \mathcal{B}_i$, with each such copy of the system being parametrized by two parameters.
\subsection{Kuramoto oscillators system}

The other system we consider is a network of Kuramoto oscillators with inertia with a tunable network. Thus while above we have a fixed adjacency matrix $A_{ij}$, here the network itself is part of the parameters and thus dentoed $p_{ij}$. The system $\mathcal{K}^{N}$ is defined as:

\begin{gather}
\ddot \phi_n = \Omega_n - p_{n} \dot \phi_n - K \sum_{m = 1, m \neq n}^N p_{nm}\sin{(x_n - x_m)} + \delta_{n1} i,\nonumber\\
n = 1 \dots N, \nonumber\\
o(t) = i(t) - \phi_1(t)\label{eq:Kuramoto system} \\ \nonumber 
\phi_n(0) = 0 \nonumber\\
p_n \in \mathbb R^+ \nonumber\\
p_{nm} \in \mathbb R^+ \nonumber\\
\nonumber 1 < n, m < N \text{ and } m \neq n
\end{gather}

\begin{gather}
\mathcal{B}_i = \left\{i : i(t) = \text{Re} \sum_{l = 0}^L a_l e^{i (2\pi l t + \theta_l)}\right\},
\end{gather}
where $K$ is the coupling constant, and $\Omega_n$ is the intrinsic node frequency drawn from a Gaussian distribution centered at $0$ with width $1$.

We have $\mathcal{S} = \mathcal{K}^{10}$ and $\mathcal{C} = \mathcal{K}^{1}$.
Thus the specification is a single oscillator.

\section{Numerical results}

We will now show the tuning for these systems. In both cases we begin by taking the system and the specification and determining the $d^\rho$ of the untuned system at randomly guessed initial parameters $p$, and then tune the system to improve $d^\rho$. We will also study the question of whether the optimization is overfitting to the sample. Overfitting here means that for a small number of samples and a large number of parameters, it can occur that the parameters fit the specific sample, rather than the underlying distribution. TO rule this out we resample after tuning the system, evaluating the final $d^\rho$ achieved on a sample different from the one used to optimize the parameters.

Specific details of tuning the two example systems are presented below. Tuning is done using the Algorithms ADAM, AMSGrad and BFGS as implemented in the package DiffEqFlux.jl. By trying different combinations of optimizers and number of iterations we found that BFGS is best suited for low dimensional parameter spaces, as occur in the estimation of $\hat d^\rho$. The tuning is a high dimensional optimization and works well with gradient descent methods, such as ADAM and AMSGrad.

Step-by-step description of the tuning algorithm and the code can be found in the project repository at \url{https://github.com/PIK-ICoNe/ProBeTune.jl}. The tuning was performed on a laptop with i7-8665U CPU. The typical runtime of the tuning is of the order of 10 minutes. 

\subsection{Diffusive non-linear network}

\begin{figure}[h]
    \centering
    \includegraphics[width=0.4\paperwidth]{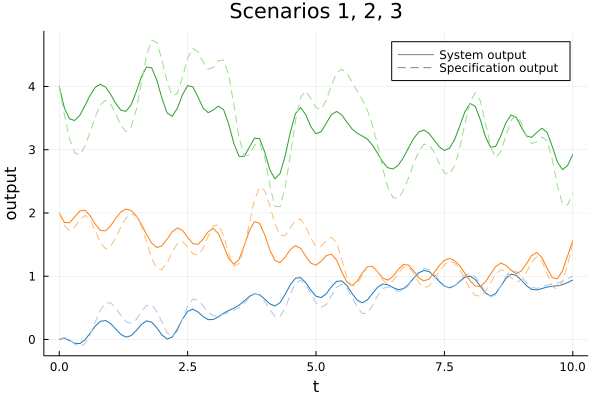}
    \caption{System and specification response to three of the inputs before tuning. $\hat d^\rho = 8.9445$. Trajectories of each element of the sample are offset by 2 on the vertical axis for better readability. Dashed line shows the response of the specification, solid line the response of the system.}
    \label{fig:init}
\end{figure}

\begin{figure}[h]
    \centering
    \includegraphics[width=0.4\paperwidth]{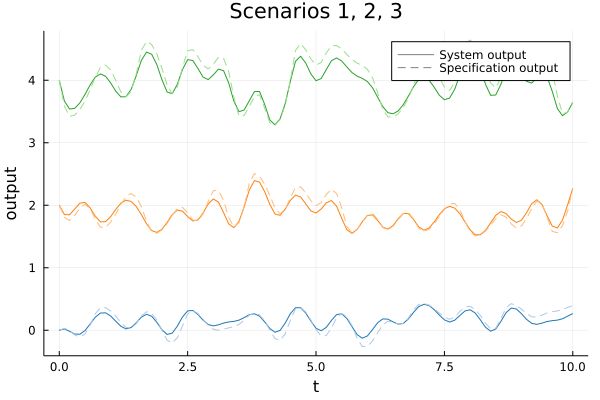}
    \caption{Results of tuning for 10 samples. $\hat d^\rho = 1.0411$. Trajectories of each element of the sample are offset by 2. Dashed line shows the response of the specification, solid line the response of the system.}
    \label{fig:10opt}
\end{figure}

Table~\ref{tab:diffusive_tuning} shows the schedule of estimating behavioral distance, resampling and tuning. We begin by estimating $\hat d^\rho$ for our initial guess of $p$ with an input sample of size ten, i.e. $|\mathcal B_i| = 10$, obtaining $\hat d^\rho \sim 8.9$. The output of system and specification for three samples is shown in Figure~\ref{fig:init}. Then we tune the system using this sample, and evaluate the distance again. This way we obtain a significant reduction in distance to the specification to $\hat d^\rho \sim 1.0$. Resampling shows that this included only moderate overfitting, the resulting outputs are shown in Figure~\ref{fig:10opt}. This might seem surprising as we only have 10 samples but should be viewed in the context that each sample consists of a timeseries and contains considerable information.

\begin{table}[h!]
\centering
\begin{tabular}{|m{7cm}|c|}
\hline
\textbf{Tuning pipeline step} & $\hat d^\rho$\\
\hline \hline
Estimate $\hat d^\rho$ with 10 samples. This provides a $q_i$ for each element of the sample. & 8.9445
\\\hline
Tuning $\hat d^\rho$ using $q_i$ from the previous step as initial parameters, 50 steps of ADAM(0.01) and 200 steps of AMSGrad(0.01) & 1.0411\\\hline

Estimating $\hat d^\rho$ on the original sample after tuning & 0.882\\\hline

Resampling the system, sample of the same size (10). Estimating $\hat d^\rho$ for the new sample. & 1.356\\\hline

Resampling the system, sample of size 100. Estimating $\hat d^\rho$ for the new sample to find initial values of the parameters. & 1.072\\\hline

10 repetitions of the following: 50 iterations of ADAM(0.01) and 200 iterations of AMSGrad(0.01) & 0.199 \\ \hline

Resampling the system, new sample of size 100. Estimating $\hat d^\rho$ for the new sample. & 0.2053\\\hline
\end{tabular}

\caption{Sequence of optimization steps in the tuning process of the diffusive non-linear system}\label{tab:diffusive_tuning}
\end{table}

\begin{figure}[h]
    \centering
    \includegraphics[width=0.4\paperwidth]{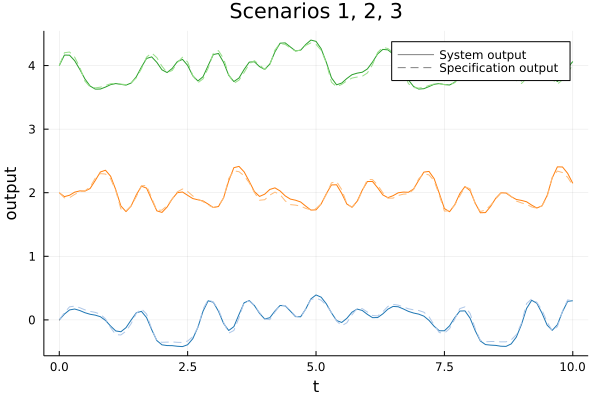}
    \caption{System and specification after tuning and resampling with 100 samples. $\hat d^\rho = 0.1978$.}
    \label{fig:complete_optim}
\end{figure}

Further tuning the system using a sample of 100 inputs shows a much further reduction of $\hat d^\rho$ to $\sim 0.2$ that persists after resampling. The outputs of the system and specification are now barely distinguishable, as shown in Figure~\ref{fig:complete_optim}. Overall the L2 distance of the difference of the output signals is a factor 45 smaller than in the untuned system.

\begin{figure}[h]
    \centering
    \includegraphics[width=.4\paperwidth]{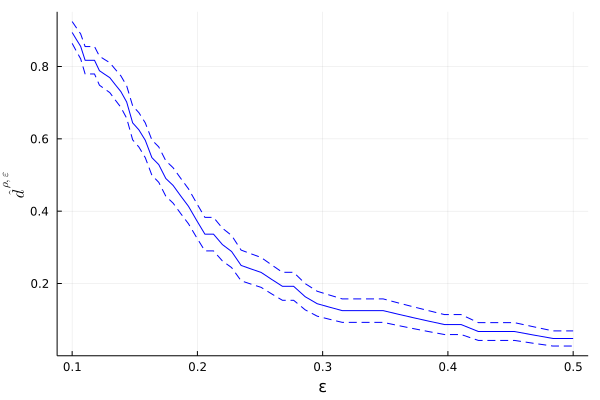}
    \caption{$\hat d^{\rho, \epsilon}$  with 95\% confidence interval}
    \label{fig:conf_int}
\end{figure}

To obtain a rigorous statement about the systems performance, we use the second distance measure \eqref{eq:estimator_distance_epsilon}. Rather than fixing a single $\epsilon$, we can plot $\hat d^{\rho, \epsilon}$ and its $95\%$ confidence interval, as a function of $\epsilon$, similar to\cite{schultz2018bounding}. This is shown in Fig.~\ref{fig:conf_int}. We thus can state that with $95\%$ confidence, more than $80\%$ of inputs produce a response that deviates less than $0.1$ from the specification in the square L2 norm chosen.

\subsection{Kuramoto oscillators system}

The Kuramoto oscillators of  Equation~\eqref{eq:Kuramoto system} are considerably more complex, exhibiting oscillating behavior. The ability to make them behave like a single oscillator will crucially depend on whether they can maintain synchrony. This in turn is controlled by the spread of the intrinsic frequencies $s$. We will see that indeed a synchronous system is easier to tune to the specification using our method.  To account for the fact that we are not tuning initial conditions to a fixed point of the system we only considered the deviation of the system and specification after initial transients have subsided.

\begin{table}[h!]
\centering
\begin{tabular}{|m{3.5cm}|c|c|}
\hline
\textbf{Tuning pipeline step} & $\hat d^\rho$ for $s=1.2$ & $\hat d^\rho$ for $s=4.5$ \\
\hline \hline
Estimate $\hat d^\rho$ with 10 samples. This provides a $q_i$ for each element of the sample. & 0.33 & 0.3853
\\\hline
Tuning $\hat d^\rho$, 100 steps of ADAM(0.01) & 0.0389 & 0.1136\\\hline
Tuning $\hat d^\rho$ further using BFGS & 0.0122 & 0.0885\\\hline
Resampling the system, sample of the same size (10). Estimating $\hat d^\rho$ for the new sample. & 0.0093 & 0.113\\\hline
Tuning $\hat d^\rho$ using BFGS & 0.0066 & -\\\hline
Resampling the system, estimating $\hat d^\rho$ for the new sample. & 0.0053 & - \\\hline
Tuning $\hat d^\rho$ further using BFGS & 0.0032 & - \\\hline
Resampling the system, estimating $\hat d^\rho$ for the new sample.& 0.0037 & - \\\hline
\end{tabular}

\caption{Sequence of optimization steps in the tuning process of the second order Kuramoto oscillators system}\label{tab:kuramoto_tuning}
\end{table}

We used a fixed set of $\Omega_n$ with $\langle \omega_n\rangle =0$ and $\omega_1=0$ and otherwise drawn from a Gaussian distribution of variance $1$, scaling them by a factor $s$. At $s = 1$ the system exhibits robust synchrony, but at $s = 5$ the system does not fully synchronize. We tuned the system with $s$ ranging from 1 to 5 with the same optimization parameters and input functions. We will show in more detail results for $s=1.2$ and $s=4.5$. Table~\ref{tab:kuramoto_tuning} provides the detailed results of two optimization schedules. For $s=4.5$ the further tuning using the BFGS optimizer failed.

\begin{figure}[h]
    \centering
    \includegraphics[width=0.4\paperwidth]{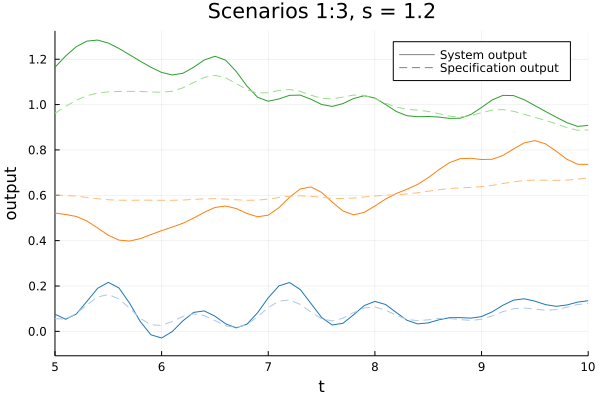}
    \caption{Initial trajectories of system and specification for $s=1.2$. $\hat d^\rho = 0.33$. We show trajectories under three inputs from the sample, offsetting them from each other by 0.5 for better readability.}
    \label{fig:kur_1.2_init}
\end{figure}

\begin{figure}[h]
    \centering
    \includegraphics[width=0.4\paperwidth]{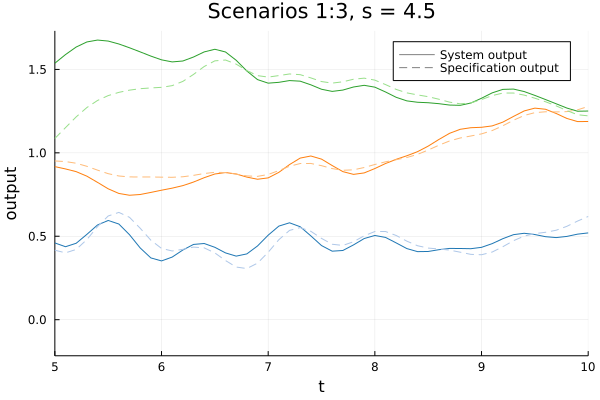}
    \caption[width=0.4\paperwidth]{Initial trajectories of system and specification for $s=4.5$. $\hat d^\rho = 0.3853$.}
    \label{fig:kur_4.5_init}
\end{figure}

Figures~\ref{fig:kur_1.2_init} and \ref{fig:kur_4.5_init} show the untuned output behavior with $\hat d^\rho$ of $\sim 0.33$ and $\sim 0.38$ respectively. 

\begin{figure}[h]
    \centering
    \includegraphics[width=0.4\paperwidth]{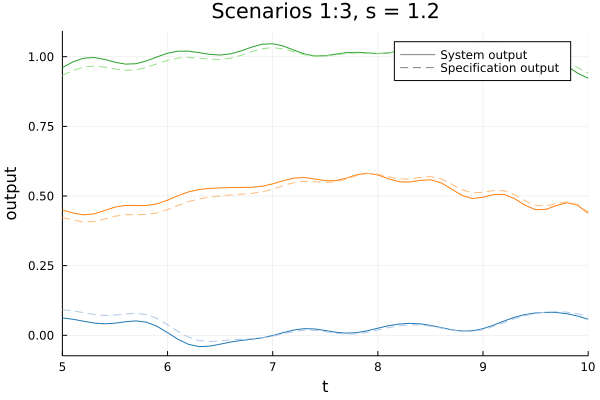}
    \caption{Results of tuning for $s=1.2$. $\hat d^\rho = 0.0032$.}
    \label{fig:kur_1.2_fin}
\end{figure}

\begin{figure}[h]
    \centering
    \includegraphics[width=0.4\paperwidth]{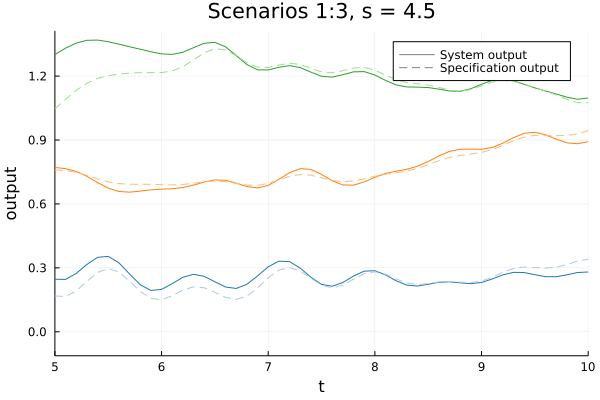}
    \caption{Results of tuning for $s=4.5$. $\hat d^\rho = 0.0885$.}
    \label{fig:kur_4.5_fin}
\end{figure}

Then we tune the systems using BFGS and ADAM algorithms. With the same schedule the respective values could be tuned down to $\sim 0.01$ and $\sim 0.11$. This shows that the synchronous behavior can be tuned much more easily in our case. Further tuning the system with $s=1.2$ we could achieve another factor three improvement of the distance, leading to an overall reduction in the square L2 norm of the outputs of a factor of 100, as compared to a factor 3 for the system with wider frequency spread. Resulting output trajectories are shown in Figure~\ref{fig:kur_1.2_fin} and Figure~\ref{fig:kur_4.5_fin}.




\begin{figure}[h]
    \centering
    \includegraphics[width=0.4\paperwidth]{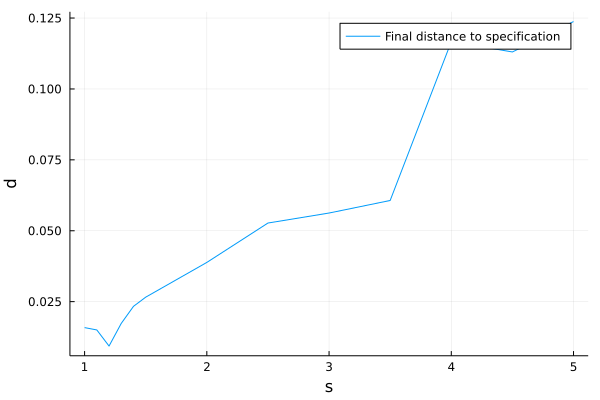}
    \caption{$\hat d^\rho$ as a function of $s$}
    \label{fig:kur_spread}
\end{figure}

To compare the effect of system synchronicity on the tuning result, we show the final value of $\hat d^\rho$ for all studied systems from s = 1 to 5 in Figure~\ref{fig:kur_spread}. For all systems we used the first four steps of the tuning pipeline of Table~\ref{tab:kuramoto_tuning} with identical parameters.

\section{Relation to $H_\infty$ model reduction}\label{sec:h-inf}

Above we defined several notions of distance between parametrized differential equations with the same inputs and outputs. These are in many ways comparable to operator norms of the difference of transfer operators of such systems. Such operator norms have been used extensively in control theory \cite{khalil1996robust}. To illuminate the similarities and differences to such approaches we will discuss the relationship of our tuning to $H_\infty$ model reduction, in which parametrized classes of transfer operators feature prominently, in more detail.

To begin with we first consider the non-probabilistic $d^{max}$ distance introduced above as this can be related explicitly to the $H_\infty$ operator norm, and its tuning to model reduction. To see this let us consider the case where both specification and system are given in terms of parametrized linear transfer operators in Laplace space $T(s)$, such that $o[i](s) = T(s)i(s)$.

The goal of $H_\infty$ model reduction is the following. Given a system with the transfer operator $T$, find a reduced system $T_{red}$ out of some class of systems, such that the difference in the induced operator norm is small. Parametrizing the reduced systems with $q$, and writing $T_{red}[q]$, we want to minimize the norm of the difference between the original and lower-order system (see e.g. \cite{chen2013robust} for an introduction): 
\begin{align}
    \min_{q} \|T-T_{red}[q]\|_\infty=\min_{q}\max_{i:\|i\|_2=1}{\|(T-T_{red}[q])i\|_2}
    \label{eq:hinf_norm}
\end{align}

We can show that the $d^{\max}$ distance of behaviors introduced in \eqref{eq:distance_max} is bounded by the $H_\infty$ norm of the optimal reduction among the $T_{red}$. Take the space of possible reduced models to be the specification $\mathcal{C}$ and the system $\overline{\mathcal S}$ given by the transfer operators $T^{\mathcal C}$ and $T^{\overline{\mathcal{S}}}$. Then we have:

\begin{align}
    d^{\max}(\overline{\mathcal{S}}, \mathcal{C}) 
    &= \max_{i : \|i\|_2 = 1} \min_q \|o^{\overline{\mathcal{S}}}[i] - o^\mathcal{C}[i, q]\|_2\nonumber\\
	&= \max_{i : \|i\|_2 = 1} \min_q {\|(T^{\overline{\mathcal{S}}} - T^{\mathcal C}[q])i\|_2}
    \label{eq:distance_hinf}
\end{align}

This differs from the norm of the optimal model reduction by the order of the $\min$ and $\max$. By the $\min$-$\max$ inequality \cite{boyd2004convex} we then have
\begin{align}
    d^{\max}(\overline{\mathcal{S}}, \mathcal{C}) 
    &\leq \min_q \max_{i : \|i\|_2 = 1} {\|(T^{\overline{\mathcal{S}}} - T^\mathcal{C}[q])i\|_2}\nonumber\\
    &= \min_q \|T^{\overline{\mathcal{S}}} - T^\mathcal{C}[q]\|_\infty
    \label{eq:hinf_inequality}
\end{align}
Thus $d^{\max}$ is a lower bound of the quality of the optimal $H_\infty$ model reduction of $T^{\overline{\mathcal{S}}}$ in the class $T^\mathcal{C}[q]$.

This comparison shows both the structural similarities and differences between standard model reduction, and our behavioral approach. Model reduction asks to have one particular simple system that behaves like the full system. Fulfilling behavioral specifications, in the setting of this paper, requires that given an input, such a system exists, but not that these simple systems are the same across all inputs. Requiring that we always have the same system is more difficult, thus the $H_\infty$ distance to the best reduction in the class under consideration is bounded by our $d^{\max}$ to the class of reductions from below.

The tuning problem we consider is to tune the complex system $\overline{\mathcal{S}}$ towards better reducability. However, both the $d^{max}$ distance and the $H_\infty$ norm require solving min-max problems that are difficult for non-linear and non-convex systems, and that might not be easy to approximate. Further, they are given by the behavior of the system given the most challenging input. Depending on the purpose of the tuning, a focus on the worst case might not be appropriate. The probabilistic distances we introduce above instead focus on the typical performance. Ignoring rare or weak failures of the specification is what allows these distances to be well approximated by sampling typical inputs. We trade a hard to track min-max problem for an easily approximated probabilistic estimation, at the price of having to provide a meaningful probability distribution on the set of inputs.

\section{Discussion}

In this paper we show how to combine probabilistic and behavioral concepts to provide novel distance measures that quantify how well a system conforms to a specification. Further, we demonstrate that these distance measures are well suited to tuning a complex system to a specification. Thus they enable us to aggregate a complex network into a vastly simpler specification. While they are probabilistic, we can give mathematically precise confidence intervals for the performance of the tuned system.

We demonstrate that the method can be efficiently used by implementing it in Julia, which has excellent library support for the type of optimizations required here. We use this implementation to successfully tune a diffusive non-linear networked system with 10 nodes to behave as a 2-node system, by jointly optimizing the system and one hundred copies of the specification over a sample of likely inputs, where each copy of specification corresponds to one possible input.

We also tune a system of 10 Kuramoto oscillators, making them behave in terms of input-output behavior as a single oscillator. This example is important for potential application of these ideas to power grids, as in that case the system is oscillating. The conceptual setup does not change, but the performance of the tuning depends on the system synchronicity. We have explored the impact of the intrinsic frequencies spread on the tuning result and find that synchronous systems are easier to tune to represent as a single oscillator, while for non-synchronous systems the same quality cannot be achieved.

The quality of tuning is evaluated using probabilistic notions. We count distances below and above a certain threshold $\varepsilon$ as successes and failures respectively, interpreting it as a Bernoulli trial to provide a rigorous confidence interval. This leads to a stability curve similar to those underlying the approach in \cite{schultz2018bounding}.

Finally we also discuss the relationship of this approach to model reduction and of our distances to $H_\infty$ norms in this context. While the use cases are different, there are considerable structural similarities.

The method introduced here is in principle well suited to establishing novel control hierarchies in complex multi-modal systems. The motivating example being future renewable power grids. In this context, the challenge is to optimize parameters of an energy cell, where the input and output characterize the power flow at the grid connection point. This sub-network should the be tuned to a specification that ensures that the continental scale system made of these cells is stable.

However, many open questions remain in order to employ this approach in such a realistic context, and the approach introduced here raises many new questions. Most importantly, in order to realize the potential of novel control hierarchies based on probabilistically satisfied specifications, we need to understand how to safely compose such specifications in a way that the guaranteed probabilistic properties are preserved.

We also focused here on specifications provided by parametrized differential equations. As noted in the introduction, the specifications for power grids usually given in terms of direct properties of the trajectory though. In this context it might be possible to explicitly solve for the specification compliant output that has the least distance to the system output. This would bypass the spec parameter optimization, and all our remaining concepts would carry through the same way.

Applying the method to non-linear systems relies on the ability of optimization algorithms to perform efficient searches in the parameter space. While the ability to differentiate through ODE solvers means that a wide variety of solvers are available for this task in Julia, the systems explored so far do not allow a comprehensive picture of their performance characteristics. Finally hand tailored optimization algorithms for this problem are also an intriguing possibility.

\section{Acknowledgments}

This research is partially funded within the framework of the "Technology-oriented systems analysis" funding area of the BMWi's 7th Energy Research Programme "Innovation for the energy transition" (FKZ: 03EI1016B) and the Deutsche Forschungsgemeinschaft (Grant No. KU 837/39-1 / RA 516/13-1).

\bibliography{main}

\providecommand{\noopsort}[1]{}\providecommand{\singleletter}[1]{#1}%
\begin{thebibliography}{23}%
\makeatletter
\providecommand \@ifxundefined [1]{%
 \@ifx{#1\undefined}
}%
\providecommand \@ifnum [1]{%
 \ifnum #1\expandafter \@firstoftwo
 \else \expandafter \@secondoftwo
 \fi
}%
\providecommand \@ifx [1]{%
 \ifx #1\expandafter \@firstoftwo
 \else \expandafter \@secondoftwo
 \fi
}%
\providecommand \natexlab [1]{#1}%
\providecommand \enquote  [1]{``#1''}%
\providecommand \bibnamefont  [1]{#1}%
\providecommand \bibfnamefont [1]{#1}%
\providecommand \citenamefont [1]{#1}%
\providecommand \href@noop [0]{\@secondoftwo}%
\providecommand \href [0]{\begingroup \@sanitize@url \@href}%
\providecommand \@href[1]{\@@startlink{#1}\@@href}%
\providecommand \@@href[1]{\endgroup#1\@@endlink}%
\providecommand \@sanitize@url [0]{\catcode `\\12\catcode `\$12\catcode
  `\&12\catcode `\#12\catcode `\^12\catcode `\_12\catcode `\%12\relax}%
\providecommand \@@startlink[1]{}%
\providecommand \@@endlink[0]{}%
\providecommand \url  [0]{\begingroup\@sanitize@url \@url }%
\providecommand \@url [1]{\endgroup\@href {#1}{\urlprefix }}%
\providecommand \urlprefix  [0]{URL }%
\providecommand \Eprint [0]{\href }%
\providecommand \doibase [0]{http://dx.doi.org/}%
\providecommand \selectlanguage [0]{\@gobble}%
\providecommand \bibinfo  [0]{\@secondoftwo}%
\providecommand \bibfield  [0]{\@secondoftwo}%
\providecommand \translation [1]{[#1]}%
\providecommand \BibitemOpen [0]{}%
\providecommand \bibitemStop [0]{}%
\providecommand \bibitemNoStop [0]{.\EOS\space}%
\providecommand \EOS [0]{\spacefactor3000\relax}%
\providecommand \BibitemShut  [1]{\csname bibitem#1\endcsname}%
\let\auto@bib@innerbib\@empty
\bibitem [{\citenamefont {Willems}(1989)}]{willems1989models}%
  \BibitemOpen
  \bibfield  {author} {\bibinfo {author} {\bibfnamefont {J.~C.}\ \bibnamefont
  {Willems}},\ }in\ \href@noop {} {\emph {\bibinfo {booktitle} {Dynamics
  reported}}}\ (\bibinfo  {publisher} {Springer},\ \bibinfo {year} {1989})\
  pp.\ \bibinfo {pages} {171--269}\BibitemShut {NoStop}%
\bibitem [{\citenamefont {Willems}\ and\ \citenamefont
  {Polderman}(1997)}]{willems1997introduction}%
  \BibitemOpen
  \bibfield  {author} {\bibinfo {author} {\bibfnamefont {J.~C.}\ \bibnamefont
  {Willems}}\ and\ \bibinfo {author} {\bibfnamefont {J.~W.}\ \bibnamefont
  {Polderman}},\ }\href@noop {} {\emph {\bibinfo {title} {Introduction to
  mathematical systems theory: a behavioral approach}}},\ Vol.~\bibinfo
  {volume} {26}\ (\bibinfo  {publisher} {Springer Science \& Business Media},\
  \bibinfo {year} {1997})\BibitemShut {NoStop}%
\bibitem [{\citenamefont {Moor}\ \emph {et~al.}(2003)\citenamefont {Moor},
  \citenamefont {Raisch},\ and\ \citenamefont
  {Davoren}}]{moor2003admissibility}%
  \BibitemOpen
  \bibfield  {author} {\bibinfo {author} {\bibfnamefont {T.}~\bibnamefont
  {Moor}}, \bibinfo {author} {\bibfnamefont {J.}~\bibnamefont {Raisch}}, \ and\
  \bibinfo {author} {\bibfnamefont {J.~M.}\ \bibnamefont {Davoren}},\
  }\href@noop {} {\bibfield  {journal} {\bibinfo  {journal} {IFAC Proceedings
  Volumes}\ }\textbf {\bibinfo {volume} {36}},\ \bibinfo {pages} {349}
  (\bibinfo {year} {2003})}\BibitemShut {NoStop}%
\bibitem [{\citenamefont {Rackauckas}\ and\ \citenamefont
  {Nie}(2017)}]{rackauckas2017differentialequations}%
  \BibitemOpen
  \bibfield  {author} {\bibinfo {author} {\bibfnamefont {C.}~\bibnamefont
  {Rackauckas}}\ and\ \bibinfo {author} {\bibfnamefont {Q.}~\bibnamefont
  {Nie}},\ }\href@noop {} {\bibfield  {journal} {\bibinfo  {journal} {Journal
  of Open Research Software}\ }\textbf {\bibinfo {volume} {5}} (\bibinfo {year}
  {2017})}\BibitemShut {NoStop}%
\bibitem [{\citenamefont {Rackauckas}\ \emph {et~al.}(2019)\citenamefont
  {Rackauckas}, \citenamefont {Innes}, \citenamefont {Ma}, \citenamefont
  {Bettencourt}, \citenamefont {White},\ and\ \citenamefont
  {Dixit}}]{rackauckas2019diffeqflux}%
  \BibitemOpen
  \bibfield  {author} {\bibinfo {author} {\bibfnamefont {C.}~\bibnamefont
  {Rackauckas}}, \bibinfo {author} {\bibfnamefont {M.}~\bibnamefont {Innes}},
  \bibinfo {author} {\bibfnamefont {Y.}~\bibnamefont {Ma}}, \bibinfo {author}
  {\bibfnamefont {J.}~\bibnamefont {Bettencourt}}, \bibinfo {author}
  {\bibfnamefont {L.}~\bibnamefont {White}}, \ and\ \bibinfo {author}
  {\bibfnamefont {V.}~\bibnamefont {Dixit}},\ }\href@noop {} {\bibfield
  {journal} {\bibinfo  {journal} {arXiv preprint arXiv:1902.02376}\ } (\bibinfo
  {year} {2019})}\BibitemShut {NoStop}%
\bibitem [{\citenamefont {B\"uttner}\ \emph {et~al.}(2021)\citenamefont
  {B\"uttner}, \citenamefont {W\"urfel}, \citenamefont {Plietzsch},
  \citenamefont {Lindner},\ and\ \citenamefont {Hellmann}}]{buttner2021stack}%
  \BibitemOpen
  \bibfield  {author} {\bibinfo {author} {\bibfnamefont {A.}~\bibnamefont
  {B\"uttner}}, \bibinfo {author} {\bibfnamefont {H.}~\bibnamefont {W\"urfel}},
  \bibinfo {author} {\bibfnamefont {A.}~\bibnamefont {Plietzsch}}, \bibinfo
  {author} {\bibfnamefont {M.}~\bibnamefont {Lindner}}, \ and\ \bibinfo
  {author} {\bibfnamefont {F.}~\bibnamefont {Hellmann}},\ }\href@noop {}
  {\bibfield  {journal} {\bibinfo  {journal} {to appear}\ } (\bibinfo {year}
  {2021})}\BibitemShut {NoStop}%
\bibitem [{\citenamefont {Lindner}\ \emph {et~al.}(2021)\citenamefont
  {Lindner}, \citenamefont {Lincoln}, \citenamefont {Drauschke}, \citenamefont
  {Koulen}, \citenamefont {Würfel}, \citenamefont {Plietzsch},\ and\
  \citenamefont {Hellmann}}]{lindner2021networkdynamics}%
  \BibitemOpen
  \bibfield  {author} {\bibinfo {author} {\bibfnamefont {M.}~\bibnamefont
  {Lindner}}, \bibinfo {author} {\bibfnamefont {L.}~\bibnamefont {Lincoln}},
  \bibinfo {author} {\bibfnamefont {F.}~\bibnamefont {Drauschke}}, \bibinfo
  {author} {\bibfnamefont {J.~M.}\ \bibnamefont {Koulen}}, \bibinfo {author}
  {\bibfnamefont {H.}~\bibnamefont {Würfel}}, \bibinfo {author} {\bibfnamefont
  {A.}~\bibnamefont {Plietzsch}}, \ and\ \bibinfo {author} {\bibfnamefont
  {F.}~\bibnamefont {Hellmann}},\ }\href {\doibase 10.1063/5.0051387}
  {\bibfield  {journal} {\bibinfo  {journal} {Chaos: An Interdisciplinary
  Journal of Nonlinear Science}\ }\textbf {\bibinfo {volume} {31}},\ \bibinfo
  {pages} {063133} (\bibinfo {year} {2021})},\ \Eprint
  {http://arxiv.org/abs/https://doi.org/10.1063/5.0051387}
  {https://doi.org/10.1063/5.0051387} \BibitemShut {NoStop}%
\bibitem [{\citenamefont {Plietzsch}\ \emph {et~al.}(2021)\citenamefont
  {Plietzsch}, \citenamefont {Kogler}, \citenamefont {Auer}, \citenamefont
  {Merino}, \citenamefont {Gil-de Muro}, \citenamefont {Li{\ss}e},
  \citenamefont {Vogel},\ and\ \citenamefont
  {Hellmann}}]{plietzsch2021powerdynamics}%
  \BibitemOpen
  \bibfield  {author} {\bibinfo {author} {\bibfnamefont {A.}~\bibnamefont
  {Plietzsch}}, \bibinfo {author} {\bibfnamefont {R.}~\bibnamefont {Kogler}},
  \bibinfo {author} {\bibfnamefont {S.}~\bibnamefont {Auer}}, \bibinfo {author}
  {\bibfnamefont {J.}~\bibnamefont {Merino}}, \bibinfo {author} {\bibfnamefont
  {A.}~\bibnamefont {Gil-de Muro}}, \bibinfo {author} {\bibfnamefont
  {J.}~\bibnamefont {Li{\ss}e}}, \bibinfo {author} {\bibfnamefont
  {C.}~\bibnamefont {Vogel}}, \ and\ \bibinfo {author} {\bibfnamefont
  {F.}~\bibnamefont {Hellmann}},\ }\href@noop {} {\bibfield  {journal}
  {\bibinfo  {journal} {arXiv preprint arXiv:2101.02103}\ } (\bibinfo {year}
  {2021})}\BibitemShut {NoStop}%
\bibitem [{\citenamefont {ENTSO-E}(2018)}]{entsoe}%
  \BibitemOpen
  \bibfield  {author} {\bibinfo {author} {\bibnamefont {ENTSO-E}},\ }\href
  {https://www.entsoe.eu/network_codes/cnc/cnc-igds/} {\enquote {\bibinfo
  {title} {Rate of change of frequency (rocof) withstand capability},}\ }
  (\bibinfo {year} {2018})\BibitemShut {NoStop}%
\bibitem [{\citenamefont {Borkowska}(1974)}]{borkowska1974probabilistic}%
  \BibitemOpen
  \bibfield  {author} {\bibinfo {author} {\bibfnamefont {B.}~\bibnamefont
  {Borkowska}},\ }\href@noop {} {\bibfield  {journal} {\bibinfo  {journal}
  {IEEE Transactions on Power Apparatus and Systems}\ ,\ \bibinfo {pages}
  {752}} (\bibinfo {year} {1974})}\BibitemShut {NoStop}%
\bibitem [{\citenamefont {Anders}(1989)}]{anders1989probability}%
  \BibitemOpen
  \bibfield  {author} {\bibinfo {author} {\bibfnamefont {G.~J.}\ \bibnamefont
  {Anders}},\ }\href@noop {} {\emph {\bibinfo {title} {Probability concepts in
  electric power systems}}}\ (\bibinfo  {publisher} {New York, NY; John Wiley
  and Sons Inc.},\ \bibinfo {year} {1989})\BibitemShut {NoStop}%
\bibitem [{\citenamefont {Menck}\ \emph {et~al.}(2014)\citenamefont {Menck},
  \citenamefont {Heitzig}, \citenamefont {Kurths},\ and\ \citenamefont
  {Schellnhuber}}]{menck2014dead}%
  \BibitemOpen
  \bibfield  {author} {\bibinfo {author} {\bibfnamefont {P.~J.}\ \bibnamefont
  {Menck}}, \bibinfo {author} {\bibfnamefont {J.}~\bibnamefont {Heitzig}},
  \bibinfo {author} {\bibfnamefont {J.}~\bibnamefont {Kurths}}, \ and\ \bibinfo
  {author} {\bibfnamefont {H.~J.}\ \bibnamefont {Schellnhuber}},\ }\href@noop
  {} {\bibfield  {journal} {\bibinfo  {journal} {Nature communications}\
  }\textbf {\bibinfo {volume} {5}},\ \bibinfo {pages} {1} (\bibinfo {year}
  {2014})}\BibitemShut {NoStop}%
\bibitem [{\citenamefont {Hellmann}\ \emph {et~al.}(2016)\citenamefont
  {Hellmann}, \citenamefont {Schultz}, \citenamefont {Grabow}, \citenamefont
  {Heitzig},\ and\ \citenamefont {Kurths}}]{hellmann2016survivability}%
  \BibitemOpen
  \bibfield  {author} {\bibinfo {author} {\bibfnamefont {F.}~\bibnamefont
  {Hellmann}}, \bibinfo {author} {\bibfnamefont {P.}~\bibnamefont {Schultz}},
  \bibinfo {author} {\bibfnamefont {C.}~\bibnamefont {Grabow}}, \bibinfo
  {author} {\bibfnamefont {J.}~\bibnamefont {Heitzig}}, \ and\ \bibinfo
  {author} {\bibfnamefont {J.}~\bibnamefont {Kurths}},\ }\href@noop {}
  {\bibfield  {journal} {\bibinfo  {journal} {Scientific reports}\ }\textbf
  {\bibinfo {volume} {6}},\ \bibinfo {pages} {1} (\bibinfo {year}
  {2016})}\BibitemShut {NoStop}%
\bibitem [{\citenamefont {Hellmann}\ \emph {et~al.}(2020)\citenamefont
  {Hellmann}, \citenamefont {Schultz}, \citenamefont {Jaros}, \citenamefont
  {Levchenko}, \citenamefont {Kapitaniak}, \citenamefont {Kurths},\ and\
  \citenamefont {Maistrenko}}]{hellmann2020network}%
  \BibitemOpen
  \bibfield  {author} {\bibinfo {author} {\bibfnamefont {F.}~\bibnamefont
  {Hellmann}}, \bibinfo {author} {\bibfnamefont {P.}~\bibnamefont {Schultz}},
  \bibinfo {author} {\bibfnamefont {P.}~\bibnamefont {Jaros}}, \bibinfo
  {author} {\bibfnamefont {R.}~\bibnamefont {Levchenko}}, \bibinfo {author}
  {\bibfnamefont {T.}~\bibnamefont {Kapitaniak}}, \bibinfo {author}
  {\bibfnamefont {J.}~\bibnamefont {Kurths}}, \ and\ \bibinfo {author}
  {\bibfnamefont {Y.}~\bibnamefont {Maistrenko}},\ }\href@noop {} {\bibfield
  {journal} {\bibinfo  {journal} {Nature communications}\ }\textbf {\bibinfo
  {volume} {11}},\ \bibinfo {pages} {1} (\bibinfo {year} {2020})}\BibitemShut
  {NoStop}%
\bibitem [{\citenamefont {Liemann}\ \emph {et~al.}(2020)\citenamefont
  {Liemann}, \citenamefont {Strenge}, \citenamefont {Schultz}, \citenamefont
  {Hinners}, \citenamefont {Porst}, \citenamefont {Sarstedt},\ and\
  \citenamefont {Hellmann}}]{liemann2020probabilistic}%
  \BibitemOpen
  \bibfield  {author} {\bibinfo {author} {\bibfnamefont {S.}~\bibnamefont
  {Liemann}}, \bibinfo {author} {\bibfnamefont {L.}~\bibnamefont {Strenge}},
  \bibinfo {author} {\bibfnamefont {P.}~\bibnamefont {Schultz}}, \bibinfo
  {author} {\bibfnamefont {H.}~\bibnamefont {Hinners}}, \bibinfo {author}
  {\bibfnamefont {J.}~\bibnamefont {Porst}}, \bibinfo {author} {\bibfnamefont
  {M.}~\bibnamefont {Sarstedt}}, \ and\ \bibinfo {author} {\bibfnamefont
  {F.}~\bibnamefont {Hellmann}},\ }\href@noop {} {\bibfield  {journal}
  {\bibinfo  {journal} {IEEE Madrid PowerTech 2021}\ } (\bibinfo {year}
  {2020})}\BibitemShut {NoStop}%
\bibitem [{\citenamefont {Filip}\ \emph {et~al.}(2019)\citenamefont {Filip},
  \citenamefont {Javeed},\ and\ \citenamefont {Trefethen}}]{filip2019smooth}%
  \BibitemOpen
  \bibfield  {author} {\bibinfo {author} {\bibfnamefont {S.}~\bibnamefont
  {Filip}}, \bibinfo {author} {\bibfnamefont {A.}~\bibnamefont {Javeed}}, \
  and\ \bibinfo {author} {\bibfnamefont {L.~N.}\ \bibnamefont {Trefethen}},\
  }\href@noop {} {\bibfield  {journal} {\bibinfo  {journal} {SIAM Review}\
  }\textbf {\bibinfo {volume} {61}},\ \bibinfo {pages} {185} (\bibinfo {year}
  {2019})}\BibitemShut {NoStop}%
\bibitem [{\citenamefont {Agresti}\ and\ \citenamefont
  {Coull}(1998)}]{agresti1998approximate}%
  \BibitemOpen
  \bibfield  {author} {\bibinfo {author} {\bibfnamefont {A.}~\bibnamefont
  {Agresti}}\ and\ \bibinfo {author} {\bibfnamefont {B.~A.}\ \bibnamefont
  {Coull}},\ }\href@noop {} {\bibfield  {journal} {\bibinfo  {journal} {The
  American Statistician}\ }\textbf {\bibinfo {volume} {52}},\ \bibinfo {pages}
  {119} (\bibinfo {year} {1998})}\BibitemShut {NoStop}%
\bibitem [{\citenamefont {Bezanson}\ \emph {et~al.}(2017)\citenamefont
  {Bezanson}, \citenamefont {Edelman}, \citenamefont {Karpinski},\ and\
  \citenamefont {Shah}}]{bezanson2017julia}%
  \BibitemOpen
  \bibfield  {author} {\bibinfo {author} {\bibfnamefont {J.}~\bibnamefont
  {Bezanson}}, \bibinfo {author} {\bibfnamefont {A.}~\bibnamefont {Edelman}},
  \bibinfo {author} {\bibfnamefont {S.}~\bibnamefont {Karpinski}}, \ and\
  \bibinfo {author} {\bibfnamefont {V.~B.}\ \bibnamefont {Shah}},\ }\href@noop
  {} {\bibfield  {journal} {\bibinfo  {journal} {SIAM review}\ }\textbf
  {\bibinfo {volume} {59}},\ \bibinfo {pages} {65} (\bibinfo {year}
  {2017})}\BibitemShut {NoStop}%
\bibitem [{\citenamefont {Barab{\'a}si}\ and\ \citenamefont
  {Albert}(1999)}]{barabasi1999emergence}%
  \BibitemOpen
  \bibfield  {author} {\bibinfo {author} {\bibfnamefont {A.-L.}\ \bibnamefont
  {Barab{\'a}si}}\ and\ \bibinfo {author} {\bibfnamefont {R.}~\bibnamefont
  {Albert}},\ }\href@noop {} {\bibfield  {journal} {\bibinfo  {journal}
  {science}\ }\textbf {\bibinfo {volume} {286}},\ \bibinfo {pages} {509}
  (\bibinfo {year} {1999})}\BibitemShut {NoStop}%
\bibitem [{\citenamefont {Schultz}\ \emph {et~al.}(2018)\citenamefont
  {Schultz}, \citenamefont {Hellmann}, \citenamefont {Webster},\ and\
  \citenamefont {Kurths}}]{schultz2018bounding}%
  \BibitemOpen
  \bibfield  {author} {\bibinfo {author} {\bibfnamefont {P.}~\bibnamefont
  {Schultz}}, \bibinfo {author} {\bibfnamefont {F.}~\bibnamefont {Hellmann}},
  \bibinfo {author} {\bibfnamefont {K.~N.}\ \bibnamefont {Webster}}, \ and\
  \bibinfo {author} {\bibfnamefont {J.}~\bibnamefont {Kurths}},\ }\href@noop {}
  {\bibfield  {journal} {\bibinfo  {journal} {Chaos: An Interdisciplinary
  Journal of Nonlinear Science}\ }\textbf {\bibinfo {volume} {28}},\ \bibinfo
  {pages} {043102} (\bibinfo {year} {2018})}\BibitemShut {NoStop}%
\bibitem [{\citenamefont {Khalil}\ \emph {et~al.}(1996)\citenamefont {Khalil},
  \citenamefont {Doyle},\ and\ \citenamefont {Glover}}]{khalil1996robust}%
  \BibitemOpen
  \bibfield  {author} {\bibinfo {author} {\bibfnamefont {I.~S.}\ \bibnamefont
  {Khalil}}, \bibinfo {author} {\bibfnamefont {J.}~\bibnamefont {Doyle}}, \
  and\ \bibinfo {author} {\bibfnamefont {K.}~\bibnamefont {Glover}},\
  }\href@noop {} {\emph {\bibinfo {title} {Robust and optimal control}}}\
  (\bibinfo  {publisher} {prentice hall, new jersey},\ \bibinfo {year}
  {1996})\BibitemShut {NoStop}%
\bibitem [{\citenamefont {Chen}(2013)}]{chen2013robust}%
  \BibitemOpen
  \bibfield  {author} {\bibinfo {author} {\bibfnamefont {B.~M.}\ \bibnamefont
  {Chen}},\ }\href@noop {} {\emph {\bibinfo {title} {Robust and H-$\infty$
  Control}}}\ (\bibinfo  {publisher} {Springer Science \& Business Media},\
  \bibinfo {year} {2013})\BibitemShut {NoStop}%
\bibitem [{\citenamefont {Boyd}\ \emph {et~al.}(2004)\citenamefont {Boyd},
  \citenamefont {Boyd},\ and\ \citenamefont {Vandenberghe}}]{boyd2004convex}%
  \BibitemOpen
  \bibfield  {author} {\bibinfo {author} {\bibfnamefont {S.}~\bibnamefont
  {Boyd}}, \bibinfo {author} {\bibfnamefont {S.~P.}\ \bibnamefont {Boyd}}, \
  and\ \bibinfo {author} {\bibfnamefont {L.}~\bibnamefont {Vandenberghe}},\
  }\href@noop {} {\emph {\bibinfo {title} {Convex optimization}}}\ (\bibinfo
  {publisher} {Cambridge university press},\ \bibinfo {year}
  {2004})\BibitemShut {NoStop}%
\end{thebibliography}%


\providecommand{\noopsort}[1]{}\providecommand{\singleletter}[1]{#1}%
%
\newpage



\end{document}